\documentclass[11pt,twoside]{article}

\usepackage{a4wide}
\usepackage{amsfonts}
\usepackage{amssymb}
\usepackage{amsmath}
\usepackage{graphicx}

\newcommand{\ignore}[1]{}

\def\@begintheorem#1#2{\par\bgroup{\sc #1\ #2. }\it\ignorespaces}
\def\@opargbegintheorem#1#2#3{\par\bgroup{\sc #1\ #2\ (#3). } \it\ignorespaces}
\def\@endtheorem{\egroup}
\newtheorem{theorem}{Theorem}[section]
\newtheorem{corollary}[theorem]{Corollary}
\newtheorem{lemma}[theorem]{Lemma}

\newtheorem{example}[theorem]{Example}
\newtheorem{proposition}[theorem]{Proposition}
\newtheorem{definition}[theorem]{Definition}
\newcommand{\bt}[1]{\begin{theorem}\label{#1}}
\newcommand{\bc}[1]{\begin{corollary}\label{#1}}
\newcommand{\bl}[1]{\begin{lemma}\label{#1}}
\newcommand{\be}[1]{\begin{example}\label{#1}}
\newcommand{\bp}[1]{\begin{proposition}\label{#1}}
\newcommand{\ba}[1]{\begin{algorithm}\rm\label{#1}}
\newcommand{\bd}[1]{\begin{definition}\rm\label{#1}}{\normalsize }

\newcommand{\et}{\end{theorem}}
\newcommand{\ec}{\end{corollary}}
\newcommand{\el}{\end{lemma}}
\newcommand{\ee}{\end{example}}
\newcommand{\ep}{\end{proposition}}
\newcommand{\ed}{\end{definition}}
\newcommand{\epr}{{\ \vbox{\hrule\hbox{%
\vrule height1.3ex\hskip0.8ex\vrule}\hrule}}\\\par}

\def\Z{\mathbb{Z}}
\def\l{\lambda}
\def\L{\Lambda}

\begin{document}

\title{\bf Optimization over Young Diagrams}

\author{
Shmuel Onn
\thanks{\small Technion - Israel Institute of Technology.
Email: onn@technion.ac.il}
}
\date{}

\maketitle

\begin{abstract}
We consider the problem of finding a Young diagram minimizing the sum of evaluations
of a given pair of functions on the parts of the associated pair of conjugate partitions.
While there are exponentially many diagrams, we show it is polynomial time solvable.

\end{abstract}

\section{Introduction}

For a Young diagram $\L$, let $n=|\L|$ be the number of cells, let $\l\vdash n$ be the partition
of $n$ whose $i$-th part $\l_i$ is the number of cells in the $i$-th row, and let $\l^*\vdash n$
be the conjugate partition of $n$ whose $j$-th part $\l^*_j$ is the number of cells in the
$j$-th column, so that $\l^*_j=|\{i:\l_i\geq j\}|$. For function $f:[n]\rightarrow\Z$ let
$f(\l)=\sum_i f(\l_i)$ be the sum of evaluations of $f$ on the parts of $\l$.
See \cite{Sta} for more information on Young diagrams and partitions and their many applications.

We consider here the following algorithmic problem.

\vskip.2cm\noindent{\bf Optimization over Young Diagrams.}
Given $n$ and functions $f,f^*:[n]\rightarrow\Z$, find a Young diagram $\L$ which
minimizes $f(\l)+f^*(\l^*)$. Equivalently, solve $\min\{f(\l)+f^*(\l^*)\,:\,\l\vdash n\}$.

\be{example}
Let $n=6$ and $f(k)=f^*(k)=k^2$. Then there are $11$ Young diagrams $\L$ with
$$\l\ =\ (6),\ (5,1),\ (4,2),\ (4,1^2),\ (3^2),\ (3,2,1),\
(3,1^3),\ (2^3),\ (2^2,1^2),\ (2,1^4),\ (1^6)\ ,$$
$$\l^*\ =\ (1^6),\ (2,1^4),\ (2^2,1^2),\ (3,1^3),\ (2^3),\
(3,2,1),\ (4,1^2),\ (3^2),\ (4,2),\ (5,1),\ (6)\ .$$
Computing the objective function $f(\l)+f^*(\l^*)$ exhaustively for all we find that
the unique optimal one is the self conjugate $\l=\l^*=(3,2,1)$ with value
$(3^2+2^2+1^2)+(3^2+2^2+1^2)=28$.
\ee

The number of Young diagrams is exponential in $n$ and so solution by exhaustive search in
general is prohibitive. Nonetheless, we show that the problem is polynomial time solvable.

\bt{main}
Optimization over Young diagrams can be done in time polynomial in $n$.
\et

\section{Proof}

We begin with a construction of Young diagrams which will be necessary for our purposes.
The {\em type} of a Young diagram $\L$ and of the associated partition $\l$ is the number of
distinct parts of $\l$. It is easy to see that $\l$ and $\l^*$ have the same type which is equal
to the number of ``southeast corners" of $\L$. If $|\L|=n$ and $\L$ has type $k$ so
$\l=(r_1^{c_1},\dots,r_k^{c_k})$ for some $r_1>\cdots>r_k\geq1$ and some
$c_1,\dots,c_k\geq1$ then $n=\sum_{i=1}^kc_ir_i\geq\sum_{i=1}^k i$ so $k<\sqrt{2n}$. For
instance, $\l=(19,\dots,2,1^{11})$ is a partition of $n=200$ of maximum type $k=19<20=\sqrt{2n}$.

Let $1\leq k<\sqrt{2n}$ and let $n\geq r_1>\cdots>r_k>r_{k+1}=0$ and $0=c_0<c_1<\cdots<c_k\leq n$.
These numbers define the Young diagram $\L$ of type $k$ which, for $i=1,\dots,k$, has
$c_i-c_{i-1}$ rows with $r_i$ cells and $r_i-r_{i+1}$ columns with $c_i$ cells,
with partition and conjugate partition
$$\l=(r_1^{c_1-c_0},\dots,r_k^{c_k-c_{k-1}})\ ,\quad
\l^*=(c_k^{r_k-r_{k+1}},\dots,c_1^{r_1-r_2})\ .$$
Note that $|\L|=\sum_{i=1}^k(c_i-c_{i-1})r_i=\sum_{i=1}^k(r_i-r_{i+1})c_i$ is not necessarily
equal to $n$, but any diagram with $|\L|=n$ does arise that way for a unique choice
of type $k$ and such $r_i$ and $c_j$.

\vskip.2cm
Let now $n$, $f$, $f^*$ be given. Fix any $1\leq k<\sqrt{2n}$.
We reduce the problem of finding a diagram $\L$ with $|\L|=n$ of type $k$ with
minimum $f(\l)+f^*(\l^*)$ to finding a shortest directed path in
a directed graph $D$ where each edge has a length. We construct $D$ as follows.

There are two vertices $s,t$, and vertices labeled by quadruples of
integers $(i,c_i,r_{i+1},n_i)$ for $0\leq i\leq k$, with
$1\leq c_i,r_i,n_i\leq n$ for $1\leq i\leq k$, $c_0=n_0=r_{k+1}=0$, and $n_k=n$.

There are edges $[s,(0,0,r_1,0)]$ for $1\leq r_1\leq n$ and edges
$[(k,c_k,0,n),t]$ for $1\leq c_k\leq n$, all of length $0$, and there are
edges $[(i-1,c_{i-1},r_i,n_{i-1}),(i,c_i,r_{i+1},n_i)]$ for $1\leq i\leq k$,
$c_i>c_{i-1}$, $r_{i+1}<r_i$, $n_i=n_{i-1}+(c_i-c_{i-1})r_i$, of length
$(c_i-c_{i-1})f(r_i)+(r_i-r_{i+1})f^*(c_i)$.

\vskip.2cm
Consider any directed path from $s$ to $t$ in $D$, which by the definition of $D$ looks like
$$s\longrightarrow(0,c_0=0,r_1,n_0=0)\longrightarrow(1,c_1,r_2,n_1)\longrightarrow\cdots\cdots
\longrightarrow(k,c_k,r_{k+1}=0,n_k=n)\longrightarrow t\ .$$
By definition of $D$ we have $0=c_0<c_1<\cdots<c_k\leq n$ and $n\geq r_1>\cdots>r_k>r_{k+1}=0$,
giving a Young diagram $\L$ of type $k$ as explained above, with partition and conjugate partition
$$\l=(r_1^{c_1-c_0},\dots,r_k^{c_k-c_{k-1}})\ ,\quad
\l^*=(c_k^{r_k-r_{k+1}},\dots,c_1^{r_1-r_2})\ .$$
Moreover, we have $|\L|=\sum_{i=1}^k(c_i-c_{i-1})r_i=n_k=n$, and the length of the path is
$$\sum_{i=1}^k(c_i-c_{i-1})f(r_i)+(r_i-r_{i+1})f^*(c_i)=f(\l)+f^*(\l^*)\ .$$
Conversely, it is clear that any diagram $\L$ with $|\L|=n$ of type $k$ gives such a
directed path, of length $f(\l)+f^*(\l^*)$.
So a shortest $s-t$ path gives an optimal Young diagram of type $k$.

\vskip.2cm
Now, the number of vertices of $D$ is bounded by $2+2n+(k-1)n^3=O(n^{3.5})$ and hence by
a polynomial in $n$. So a shortest directed path from $s$ to $t$ in $D$ can be obtained
in polynomial time, see e.g. \cite{Sch}, by the following simple algorithm.
For $i=1,\dots,k$, compute for every vertex $v=(i,c_i,r_{i+1},n_i)$ the length of
a shortest $s-v$ path and an edge entering it on such a shortest path,
using the values already computed for $i-1$, and then do the same for $t$.

\vskip.2cm
Now repeat the above procedure for $k=1,\dots,\lfloor\sqrt{2n}\rfloor$
and output the best diagram.
\epr

\vskip.2cm
As a simple example, for $n=2$ and $k=1$, the directed graph has exactly two $s-t$ paths,
$$s\stackrel{0}{\longrightarrow}(0,0,1,0)
\stackrel{2f(1)+f^*(2)}{\longrightarrow}(1,2,0,2)\stackrel{0}{\longrightarrow}t\ ,
\quad\mbox{of length}\quad 2f(1)+f^*(2)\ ,$$
$$s\stackrel{0}{\longrightarrow}(0,0,2,0)
\stackrel{f(2)+2f^*(1)}{\longrightarrow}(1,1,0,2)\stackrel{0}{\longrightarrow}t\ ,
\quad\mbox{of length}\quad f(2)+2f^*(1)\ ,$$
which correspond to the two diagrams $\L$ and their partitions $\l=(1^2)$ and $\l=(2)$ respectively.

\vskip.2cm
We note that the above proof in fact shows that we can solve in polynomial time the more
refined problem where we search for a best diagram among those of a prescribed type $k$.

\section*{Acknowledgments}

S. Onn was supported by a grant from the Israel Science Foundation and the Dresner chair.

\end{document}